\documentclass{IEEEtran4PSCC}

\usepackage[numbers,square]{natbib}
\ifCLASSINFOpdf
   \usepackage[pdftex]{graphicx}
\else
   \usepackage[dvips]{graphicx}
\fi
\usepackage[cmex10]{amsmath}
\usepackage{bm}
\usepackage{amssymb,amsthm}

\usepackage{tabularx}
\usepackage{booktabs}
\usepackage[table]{xcolor}

\usepackage{accents}
\usepackage{xfrac}
\usepackage{url}
\usepackage{multirow}

\usepackage{tikz}
\usetikzlibrary{calc}
\usetikzlibrary{patterns, patterns.meta}
\usetikzlibrary{arrows, arrows.meta}
\usetikzlibrary{shapes.geometric}
\usetikzlibrary {positioning}
\tikzset{vertex/.style = {shape=circle,draw,minimum size=2.7em}}
\tikzset{edge/.style = {->,> = latex'}}
\definecolor {processblue}{cmyk}{0.96,0,0,0}
\usetikzlibrary{shadows}
\usetikzlibrary{matrix, calc, arrows}
\usetikzlibrary{decorations.pathreplacing
}

\usetikzlibrary{arrows}
\usetikzlibrary{shapes.geometric}
\usetikzlibrary {positioning}
\usetikzlibrary{arrows,shapes,positioning,shadows,trees}
\tikzset{vertex/.style = {shape=circle,draw,minimum size=2.7em}}
\tikzset{edge/.style = {->,> = latex'}}
\definecolor {processblue}{cmyk}{0.96,0,0,0}
\usetikzlibrary{shadows}
\usetikzlibrary{matrix, calc, arrows}
\usetikzlibrary{decorations.pathreplacing}
\usetikzlibrary{patterns}
\usetikzlibrary{arrows.meta, chains, fit, shapes}

\newcommand\ppbb{path picture bounding box} 

\tikzset{
  basic/.style  = {every node/.style=draw, rectangle, align=center, font=\footnotesize, level distance=1.5cm},
  pruned/.style  = {every node/.style=draw, rectangle, align=center, font=\footnotesize},
  level 0/.style = {basic, rounded corners=0pt, thin,   align=center, fill=black!0},
  level 1/.style = {basic, rounded corners=0pt, thin,   align=center, fill=black!0, sibling distance = 120mm},
  level 2/.style = {basic,                      thin,   align=center, fill=black!0, sibling distance = 60mm},
  level 3/.style = {basic,                      thin,   align=center, fill=black!0, sibling distance = 35mm},
  level 4/.style = {basic,                      thin,   align=center, fill=black!0, sibling distance =31mm},
  co/.style = {
    path picture={\draw[opacity=0.3] (\ppbb.north west) -- (\ppbb.south east)
    (\ppbb.south west) -- (\ppbb.north east);}},
  pruned/.style = {opacity=1, fill=black!0, co},
  >=latex
}

\usepackage[noend,ruled,vlined]{algorithm2e}
\usepackage[hang,flushmargin]{footmisc}
\makeatletter
\newcommand{\algorithmfootnote}[2][\scriptsize]{%
  \let\old@algocf@finish\@algocf@finish 
  \def\@algocf@finish{\old@algocf@finish 
    \vspace{-2.5pt}
    \leavevmode\rlap{
    \begin{minipage}{\textwidth}
        #1#2
    \end{minipage}}%
  }%
}

\makeatletter
\let\old@ps@headings\ps@headings
\let\old@ps@IEEEtitlepagestyle\ps@IEEEtitlepagestyle
\def\psccfooter#1{%
    \def\ps@headings{%
        \old@ps@headings%
        \def\@oddfoot{\strut\hfill#1\hfill\strut}%
        \def\@evenfoot{\strut\hfill#1\hfill\strut}%
    }%
    \def\ps@IEEEtitlepagestyle{%
        \old@ps@IEEEtitlepagestyle%
        \def\@oddfoot{\strut\hfill#1\hfill\strut}%
        \def\@evenfoot{\strut\hfill#1\hfill\strut}%
    }%
    \ps@headings%
}
\makeatother

\psccfooter{%
}

\begin{document}

\newcommand{\todo}[1]{{\leavevmode{\color{blue}{\texttt{TODO} #1}}}}
\newcommand*{\red}{\textcolor{black}}
\newcommand*{\blue}{\textcolor{black}}
\title{An Efficient Hybrid Heuristic for the Transmission Expansion Planning under Uncertainty} 

\author{\IEEEauthorblockN{Yure Rocha\IEEEauthorrefmark{1},
Teobaldo Bulhões\IEEEauthorrefmark{2},
Anand Subramanian\IEEEauthorrefmark{3} and
Joaquim Dias Garcia\IEEEauthorrefmark{4}}
\IEEEauthorblockA{\IEEEauthorrefmark{1} Programa de Pós-Graduação em Informática \\
Universidade Federal da Paraíba, João Pessoa, Brazil 
}
\IEEEauthorblockA{\IEEEauthorrefmark{2} Departamento de Sistemas de Computação \\
Universidade Federal da Paraíba, João Pessoa, Brazil 
}
\IEEEauthorblockA{\IEEEauthorrefmark{3} Departamento de Computação Científica \\
Universidade Federal da Paraíba, João Pessoa, Brazil 
}
\IEEEauthorblockA{\IEEEauthorrefmark{4}
Research and Development \\
PSR - Energy Consulting and Analytics, Rio de Janeiro, Brazil 
}
}

\maketitle

\begin{abstract}
We address the stochastic transmission expansion planning (STEP) problem \blue{under} uncertaint\blue{y} in renewable generation capacity and demand. STEP's objective is to minimize total \blue{transmission} investment and generation costs. To tackle the computational challenges \blue{posed by} large-scale systems, we propose a heuristic \blue{strategy} that combines the progressive hedging (PH) algorithm for scenario-wise decomposition with an integrated \blue{approach} for solving the resulting subproblems. The latter combines a destroy-and-repair operator, a beam search procedure, and a mixed-integer programming \blue{solver}. The proposed framework is evaluated on large-scale systems from the literature \blue{with} up to 10000 nodes, adapted to \blue{stochastic} scenarios \blue{using} parameters from the California test system (CATS). Compared with a non-trivial baseline algorithm that includes the \blue{same integrated approach}, the proposed PH-based \blue{method} consistently improved solution quality for the six systems considered (including CATS), achieving an average \blue{cost improvement} of \blue{5.28}\% within a 2-hour time limit.
\end{abstract}


\thanksto{\noindent This research was partly supported by the Brazilian research agency CNPq [grants 309580/2021-8, and 306501/2025-2], and by PSR.}

\vspace{-10pt}

\section*{Notation}

{
    \footnotesize
    \setlength{\tabcolsep}{2pt}
    \noindent
    \begin{tabularx}{0.49\textwidth}{>{\raggedright\arraybackslash}p{0.35cm}>{}X}
        \textbf{Sets} & \\
        $\mathcal{S}$ & Set containing the scenarios $s$ \\
        $\mathcal{B}$ & Set of buses $b$ \\
        $\mathcal{L}$ & Set of lines (a.k.a. circuits); $\mathcal{L} = \mathcal{J} \cup \mathcal{K}$ \\
        $\mathcal{J}$ & Set of existing lines $j$; $\mathcal{J}=\{(b, b', l) \mid b, b' \in \mathcal{B}, l \in \{1, \dots, J_{bb'}\}\}$ \\
        $\mathcal{J}_b$ & Set of existing lines involving bus $b$
        \\
        $\mathcal{K}$ & Set of candidate lines $k$; $\mathcal{K}=\{(b, b', l) \mid b, b' \in \mathcal{B}, l \in \{1, \dots, K_{bb'}\}\}$ \\
        $\mathcal{K}_b$ & Set of candidate lines involving bus $b$
        \\
        $\mathcal{K}'$ & Restricted set of candidate lines; $\mathcal{K'}\subset\mathcal{K}$ \\
        $\mathcal{E}$ & Set of generators $i$ \\
        $\mathcal{E}_b$ & Set of generators connected to bus $b \in \mathcal{B}$
    \end{tabularx}
    \begin{tabularx}{.49\textwidth}{>{\raggedright\arraybackslash}p{1.2cm}>{}X}
        \textbf{Data} & \\
        $\underaccent{\bar}{G}_{is}$, $\bar{G}_{is}$ & Minimum and maximum active generation capacity for generator $i \in \mathcal{E}$, respectively, at scenario $s \in \mathcal{S}$\\
        $\bar{F}\,^0_j$, $\bar{F}\,^1_k$ & Capacity of lines  $j \in \mathcal{J}, k \in \mathcal{K}$, respectively \\
        $J_{bb'}, K_{bb'}$ & Number of existing and candidate lines, respectively, between nodes $b, b' \in \mathcal{B}$ \\
        $B^0_j, B^1_k$ & Susceptances of lines $j \in \mathcal{J}, k \in \mathcal{K}$, respectively
    \end{tabularx}
    \begin{tabularx}{.49\textwidth}{>{\raggedright\arraybackslash}p{0.65cm}>{}X}
        $C^{ope}_i$ & Operational cost for generator $i \in \mathcal{E}$ \\
        $D_{bs}$ & Active demand at bus $b \in \mathcal{B}$ at scenario $s \in \mathcal{S}$ \\
        $C^{inv}_k$ & Investment cost for line $k \in \mathcal{K}$ \\
        $p_{s}$ & Probability of scenario $s \in \mathcal{S}$ \\
        $\lambda$ & Penalty factor for power flow slacks in the linear programming model \\ 
        $\beta$ & Progressive hedging threshold value for solutions \\
        $tl$ & Time limit for subproblems
    \end{tabularx}
    \begin{tabularx}{.49\textwidth}{>{\raggedright\arraybackslash}p{0.65cm}>{}X}
        $it_{dr}$ & Destroy-and-repair maximum number of iterations \\
        $it_{bs}$, $\eta, \omega$, $N, \gamma$ & Beam search parameters: maximum number of levels without improvement, multiplier for the sizes of the subsets, maximum number of branches per node, number of nodes for branching, and node branching selection slack, respectively
    \end{tabularx}
    \begin{tabularx}{.49\textwidth}{>{\raggedright\arraybackslash}p{1.2cm}>{}X}
        \multicolumn{2}{l}{\textbf{Decision variables}} \\
        $x_k$ & Binary decision variable for the investment decision in line $k \in \mathcal{K}$; $x_k$ is equal to 1 if line $k \in \mathcal{K}$ is built, and equal to 0 otherwise \\
        $f^0_{js}$, $f^1_{ks}$ & Continuous decision variables for the power flows of lines $j \in \mathcal{J}$, $k \in \mathcal{K}$, respectively, at scenario $s \in \mathcal{S}$ \\
        $y^0_{js}$, $y^1_{ks}$ & Continuous decision variables for the power flow slacks of lines $j \in \mathcal{J}$, $k \in \mathcal{K}$, respectively, at scenario $s \in \mathcal{S}$ \\ 
        $g_{is}$ & Continuous decision variable for the power flow production by generator $i \in \mathcal{E}$ at scenario $s \in \mathcal{S}$; $\underaccent{\bar}{G}_{is} \leq g_{is} \leq \bar{G}_{is}$ \\
        $g^{bus}_{bs}$ & Continuous decision variable for the power flow injection at bus $b \in \mathcal{B}$ at scenario $s \in \mathcal{S}$; $g^{bus}_{bs} = \sum_{i \in \mathcal{E}_b}g_{is}$ \\
        $\theta_{bs}$ & Continuous decision variable for the voltage phase angle at bus $b \in \mathcal{B}$ and scenario $s \in \mathcal{S}$ \\
        $\delta_{js}$ & Difference between phase angles at scenario $s \in \mathcal{S}$ for line $j=(b, b', l) \in \mathcal{L}$:  $\delta_{js} = \theta_{bs} - \theta_{b's}$
    \end{tabularx}
}

\section{Introduction}

Transmission expansion planning (TEP) involves large-scale strategic decisions about the construction or reinforcement of transmission lines to meet the growing demand while ensuring system reliability and robustness. The objective of TEP is to minimize \blue{transmission} investment and generation costs subject to power flow constraints. In this sense, TEP has considerable economic, environmental and social impacts in both the short and long term.

\red{TEP has been shown to be \(\mathcal{NP}\)-hard \citep{moulin2010transmission}.}
Since the seminal formulation proposed in \cite{garver1970transmission}, the problem has evolved to consider more realistic aspects such as uncertainties in both short- and long-term planning (e.g., regarding renewable generation and variations in production and demand). These uncertainties can be incorporated into TEP by different approaches such as multi-scenario stochastic optimization and robust optimization. The latter approach, however, may lead to overpessimistic solutions with large investment costs \cite{bagheri2016data}.

The combinatorial nature of the problem combined with multiple scenarios in the stochastic TEP (STEP) variant results in a computationally challenging problem. In fact, even generating feasible initial solutions for single-scenario TEP problems can be time-consuming for conventional mixed integer programming (MIP) solvers. To overcome these challenges, we adopt a hybrid procedure that combines an integrated approach with a decomposition-based algorithm. First, we adopt the progressive hedging (PH) algorithm \citep{rockafellar1991scenarios}, which decomposes the stochastic problem by scenario. Second, we apply heuristic methods for the subproblems associated with the scenarios. The motivation for employing heuristic methods on subproblems comes from their computational speedups reported for different classes of problems in comparison with exact approaches.

The main contributions of this paper are summarized as follows.
\begin{itemize}
    \item We devise a simple and efficient destroy-and-repair (D\&R) operator based on binary search to, starting from an initial solution with candidate lines inserted, quickly obtain high-quality feasible solutions for single-scenario TEP problems.
    \item We introduce an iterative beam search (BS) \citep{OwMorton1988} heuristic method to enhance D\&R's feasible solutions by concurrently exploring multiple promising solutions within the solution space. An integrated approach  \blue{(IA)} then combines BS, D\&R, and a mathematical solver.
    \item We propose a hybrid solution approach consisting of the PH procedure to decompose large-scale stochastic problems into single-scenario subproblems, and the \blue{IA heuristic} to address the subproblems.
    \item We report the performance of the proposed solution approach on several large-scale STEP systems ranging from a few thousand up to 10000 buses. The results demonstrate the effectiveness of our method in systematically finding improved solutions in comparison with \blue{both} a baseline \blue{algorithm} (BA) \blue{and an alternative PH procedure that addresses subproblems with a MIP solver}. In particular, an average \blue{cost improvement} of \blue{$5.28\%$} was obtained under 2 hours of time limit.
\end{itemize}

The remainder of this paper is organized as follows. Section \ref{sec:related_work} reviews a few closely related works. Section \ref{sec:math_formulation} discusses the problem and mathematical \blue{formulations}. Section \ref{sec:sol_approach} describes our solution approach. Section \ref{sec:comp_experiments} presents the computational experiments. Finally, Section \ref{sec:conclusion} concludes \blue{and discusses future research directions}.

\section{Related Work}\label{sec:related_work}

The spatial resolution of power flow networks is usually complex, with thousands of buses scattered across hundreds of zones in a country. \blue{Furthermore}, with the increasing penetration of renewable energy, the difficulty to create representative scenarios while maintaining system tractability is usually a challenging task \citep{zuluaga2024parallel}. In this context, solution approaches capable of dealing with large-\blue{scale} multi-scenario networks can build more realistic plans with an improved compromise between these two dimensions. 

PH is a scenario-based decomposition technique often applied to solve stochastic energy capacity expansion problems \citep{munoz2015scalable, soares2022integrated, zuluaga2024parallel}. PH decomposes the problem scenario-wise, allows for each subproblem to be independently solved and combines the results at each iteration. While \blue{PH is} an exact algorithm for stochastic linear problems, \blue{it can be regarded as} a heuristic method for \blue{more general} stochastic mixed-integer models \blue{\citep{lokketangen1996progressive}}.

\citet{munoz2015scalable} addressed a stochastic planning problem that optimizes transmission and generation investments. The authors employed PH to solve a model with a set o representative scenarios obtained via a scenario reduction scheme and run experiments on a 240-bus system. \citet{zuluaga2024parallel} went beyond the previous work and investigated a large-scale capacity expansion problem involving generation, transmission and storage investment. The two-stage model adopted optimizes investment decisions in the first stage and solves a multi-period optimal power flow (OPF) problem incorporated with transmission losses in the second stage. A tailored implementation of the PH algorithm, available in the \textit{mpi-sppy} Python package, was used to solve the full model, and experiments were conducted on up to 360 scenarios of an adapted version of CATS \citep{taylor2024cats}.

Heuristic methods have also been proposed for TEP problems. Examples include greedy randomized adaptive search procedures \citep{binato2002reactive}, genetic algorithms \citep{cadini2010optimal, silvasousa2011combined}, sequential approximation approaches \citep{park2013transmission}, and constructive procedures \citep{oliveira2021constructive}. As a common practice in the literature, most of these papers report tests with a single system. Moreover, the experiments \blue{were} conducted \blue{on} small-\blue{scale} instances with a few dozens of buses. For a comprehensive review on TEP, we refer the interested reader to \cite{lumbreras2016, mahdavi2019}.

Given the above, although experiments with large-scale systems exist (see, for instance, \cite{rosemberg2024learning}), we are, to the best of our knowledge, the first to report results of a combined hybrid heuristic and decomposition approach on multiple large-scale TEP instances under uncertainty.

\section{Mathematical formulation}\label{sec:math_formulation}

This section describes the STEP problem and \blue{presents} a linear programming (LP) formulation, given in \eqref{lp_model:obj}–\eqref{lp_model:end_domains}, for scenario $s \in \mathcal{S}$.  We then present the complete extensive form of the stochastic MIP model in \eqref{model:obj}–\eqref{model:x_domain}. Despite having nonlinear constraints involving absolute values, \blue{the} LP and MIP formulations can be easily linearized.

The heuristic procedures proposed in this paper control the restricted set of candidate transmission lines to build $\mathcal{K}'$ for scenario $s \in \mathcal{S}$. In particular, the $D\&R$ method exploits information about power flows at lines to guide its search for improved solutions. The LP model allows one to quickly evaluate if the candidate transmission lines in $\mathcal{K}'$ for scenario $s \in \mathcal{S}$ yield a feasible solution to the corresponding scenario.

\blue{Uncertainties are modeled using a multi-scenario approach.} Each scenario represents a possible realization of uncertainties in renewable generation capacit\blue{y} and demand levels, and is associated with probability $p_s, s \in \mathcal{S}$. The objective of STEP is to minimize total investment costs for new transmission lines and the expected generation costs, weighted by scenario probabilities. \blue{We also adopt a static model in which reinforcement decisions apply to the entire planning horizon, unlike multi-stage or dynamic frameworks.} We assume that total generation capacity is sufficient to meet the demand, i.e., the generation expansion was performed in a previous stage. To maintain tractability, we adopt the DC-OPF approximation \citep{stott2009dc}, which enables the formulation of TEP single-scenario problems with linear constraints.

In the STEP formulation, each bus can be connected to multiple generators and is associated with at most one demand. Each existing line may have one or more corresponding candidate circuit. Power flow on a line, constrained by its maximum capacity, is determined by Kirchhoff's voltage law (KVL), i.e., by the difference between the phase angles in its extreme nodes times the line susceptance. Flow conservation at each node is enforced with Kirchhoff's current law (KCL). The LP model for subproblem $s \in \mathcal{S}$ is formally defined in \eqref{lp_model:obj}--\eqref{lp_model:end_domains}.

\renewcommand{\algorithmcfname}{Model}
\LinesNotNumbered
\setlength{\algomargin}{1.5em}
\begin{algorithm}[htbp]
    \caption{The LP model for subproblem $s \in \mathcal{S}$}
    \begin{subequations}
        \vspace{-0.6em}
        \scriptsize
        \begin{align}
            \min & \sum_{i \in \mathcal{E}}C^{ope}_ig_{is} + \lambda\sum_{j \in \mathcal{J}}y^0_{js} + \lambda\sum_{k \in \mathcal{K}}y^1_{ks} & \label{lp_model:obj}
        \end{align}
        \vspace{-2.5em}
        \begin{align}
            \text{s.t.} & \nonumber \\
            & \sum_{j \in \mathcal{J}_b}f^0_{js} + \sum_{k \in \mathcal{K}_b}f^1_{ks} + g^{bus}_{bs} = D_{bs}, & b \in \mathcal{B} \label{lp_model:fk_law} \\
            & f^0_{js} = B^0_j\delta_{js}, & j \in \mathcal{J} \label{lp_model:j_kvl} \\
            & f^1_{ks} = B^1_k\delta_{ks}, & k \in \mathcal{K}' \label{lp_model:k_kvl} \\
            & |f^0_{js}| \leq y^0_{js} + \bar{F}\,^0_j, & j \in \mathcal{J} \label{lp_model:j_fub} \\
            & |f^1_{ks}| \leq y^1_{ks} + \bar{F}\,^1_k, & k \in \mathcal{K} \label{lp_model:k_fub} \\
             & \delta_{js} = \theta_{bs} - \theta_{b's}, &j=(b, b', l) \in \mathcal{L}\label{lp_model:delta_theta} \\
            & g^{bus}_{bs} = \sum_{i \in \mathcal{E}_b}g_{is}, & b \in \mathcal{B} \label{lp_model:g_bus} \\
            & \underaccent{\bar}{G}_{is} \leq g_{is} \leq \bar{G}_{is}, & i \in \mathcal{E} \label{lp_model:g_bus_domain} \\
            & \theta_{bs} \in \mathbb{R}, & b \in \mathcal{B} \label{lp_model:theta_domain} \\
            & f^0_{js} \in \mathbb{R}, & j \in \mathcal{J} \label{lp_model:f0_domain} \\
            & f^1_{ks} \in \mathbb{R}, & k \in \mathcal{K} \label{lp_model:f1_domain} \\
            & y^0_{js} \in \mathbb{R}_{\ge 0}, & j \in \mathcal{J}  \\
            & y^1_{ks} \in \mathbb{R}_{\ge 0}, & k \in \mathcal{K} \label{lp_model:end_domains}
        \end{align}
        \vspace{-2em}
    \end{subequations}
\end{algorithm}

The LP model minimizes generation costs and the penalized slack variables \eqref{lp_model:obj}. KCL is employed in constraints \eqref{lp_model:fk_law} to ensure power flow conservation at each node of the graph. Next, KVL is applied in constraints \eqref{lp_model:j_kvl} and \eqref{lp_model:k_kvl} to model power flow for existing transmission lines and candidate transmission lines in the restricted set, respectively. Constraints \eqref{lp_model:j_fub} and \eqref{lp_model:k_fub} set the capacities of the circuits. Differences in phase angles are computed in \eqref{lp_model:delta_theta} and the total power generation at each bus is enforced by constraints \eqref{lp_model:g_bus}. Finally, constraints \eqref{lp_model:g_bus_domain}--\eqref{lp_model:end_domains} model the domain of the variables. This formulation allows power flows to exceed line capacities, subject to a penalty in the objective function \eqref{lp_model:obj}, where $\lambda$ is a sufficiently large penalty factor.

The complete extensive-form MIP model for STEP is presented in \eqref{model:obj}--\eqref{model:x_domain}.
\renewcommand{\algorithmcfname}{Model}
\LinesNotNumbered
\begin{algorithm}[htbp]
    \centering
    \caption{Extensive-form MIP model}
    \begin{subequations}
        \scriptsize
        \vspace{-0.6em}
        \begin{align}
            \min & \sum_{k \in \mathcal{K}}C^{inv}_kx_k + \sum_{i \in \mathcal{E}}C^{ope}_i\sum_{s\in\mathcal{S}}p_sg_{is} 
            & \label{model:obj}
            \end{align}
            \vspace{-2.5em}
            \begin{align}
            & \text{s.t.} & \nonumber \\
            & \eqref{lp_model:fk_law}, \eqref{lp_model:j_kvl}, \eqref{lp_model:delta_theta}, \eqref{lp_model:g_bus}, \eqref{lp_model:g_bus_domain}, \eqref{lp_model:theta_domain}, \eqref{lp_model:f0_domain}, \text{and }\eqref{lp_model:f1_domain}, & s \in \mathcal{S} \\
            & |f^1_{ks} - B^1_k\delta_{ks}| \leq M_k(1 - x_k), & k \in \mathcal{K}, s \in \mathcal{S} \label{model:k_ohm_law} \\
            & |f^0_{js}| \leq \bar{F}\,^0_j, & j \in \mathcal{J}, s \in \mathcal{S} \label{model:j_fub} \\ 
            & |f^1_{ks}| \leq \bar{F}\,^1_kx_k, & k \in \mathcal{K}, s \in \mathcal{S} \label{model:k_fub} \\
            & x_k \in \{0, 1\}, & k \in \mathcal{K}\label{model:x_domain}
        \end{align}
        \vspace{-2em}
    \end{subequations}
\end{algorithm}

The objective function \eqref{model:obj} minimizes the total investment cost of new transmission lines and the expected operational cost of generation, weighted by the probabilities of the scenarios. For scenario $s \in \mathcal{S}$ in the MIP model, KCL, KVL for existing lines, the difference between voltage angles, variable $g^{bus}_{bs}$, and the domain \blue{of} variables $g_{is}$, $\theta_{bs}$, $f^0_{js}$, and $f^1_{ks}$, are modeled in a similar fashion as in the LP model. Constraints \eqref{model:k_ohm_law} model power flow for candidate transmission lines, where $M_k, k \in \mathcal{K}$, are sufficiently large numbers. Constraints \eqref{model:j_fub} and \eqref{model:k_fub} set the capacities of the lines. Finally, constraints \eqref{model:x_domain} model the domain of investment variable $x_k$, that indicates which candidate lines to construct for all scenarios. \blue{In the context of stochastic optimization, $x_k$ denotes the first-stage (here-and-now) decisions, whereas the remaining variables correspond to the second-stage (wait-and-see) decisions.}

A solution to the LP model with $\sum_{j \in \mathcal{J}}y^0_{js} + \sum_{k \in \mathcal{K}}y^1_{ks} = 0$ is feasible for scenario $s$ of the MIP model. 

\section{Solution approach}\label{sec:sol_approach}

We employ the PH algorithm to decompose the large-scale STEP problem into multiple single-scenario subproblems. Each subproblem is then solved by means of an integrated framework consisting of construction and improvement followed by an exact approach. In the construction stage, we quickly build a feasible solution based on residual flows of existing transmission lines by means of the D\&R procedure. The improvement phase applies the BS heuristic to iteratively evaluate the removal of candidate lines in a set of promising solutions. D\&R controls a list of candidate lines constructed $\mathcal{K}'$ and uses it to evaluate violations by solving the LP model of Section \ref{sec:math_formulation}. BS works in a similar fashion, but stores multiple lists of constructed candidates lines instead.

A time limit $tl$ is set for solving each subproblem. After the execution of the D\&R-BS heuristic, the solver run\blue{s} until the time limit is reached. 

\subsection{Progressive hedging algorithm}

In our PH procedure, whose main differences from conventional implementations are discussed within this section, we adopt the cost proportional approach described \blue{in} \cite{watson2011progressive}. \blue{Figure \ref{fig:ph_diagram} illustrates the main steps executed in the PH algorithm. 
\begin{figure}
    \centering
    \blue{
    \begin{tikzpicture}[
    font=\scriptsize\sffamily,
    node distance=1.4cm,
    every node/.style={draw, rectangle, rounded corners=2pt, align=center, minimum width=2.1cm, minimum height=1cm, inner sep=3pt},
    arrow/.style={->, -{Latex[length=1mm,width=1mm]}, thin}
]
    
    \node (step) {STEP system\\(multi-scenario)};
    \node (init_cost) [below of=step] {Compute \\``build-everything'' \\solution cost};
    
    \node (bs) [below=0.8cm of init_cost] {BS \\Using LP model};
    \node (dnr) [left=0.4cm of bs] {D\&R\\Using LP model};
    \node (mip_tep) [right=0.4cm of bs] {MIP-TEP};
    
    \node (gen_sols) [below=0.8cm and 0cm of bs] {Generate \emph{union}\\and \\$\beta$-\emph{intersection}};
    \node (update_ph_data) [right=0.4cm of gen_sols] {Update PH data};
    \node (eval_sols) [left=0.4cm of gen_sols] {Evaluate$^*$ \emph{union} \\and \\$\beta$-\emph{intersection}};

    \node (select) [below of=gen_sols] {Update \\incumbent};
    \node (check) [right=0.4cm of select] {Check termina-\\tion criteria};
    \node (repair) [left=0.4cm of select] {\texttt{repair}$^*$ if \\violations occur};

    \node[
        draw,
        dashed,
        rounded corners,
        inner sep=6pt,
        fit=(dnr) (bs) (mip_tep),
    ] (ia_box) {};

    
    \draw[arrow] (step) -- (init_cost);
    \draw[arrow] (init_cost) -- ($(init_cost)!0.5!(bs)$) -| (dnr);
    
    \draw[arrow] (dnr) -- (bs);
    \draw[arrow] (bs) -- (mip_tep);
    \draw[arrow] (mip_tep) -- (update_ph_data);
    \draw[arrow] (update_ph_data) -- (gen_sols);
    \draw[arrow] (gen_sols) -- (eval_sols);
    \draw[arrow] (eval_sols) -- (repair);
    \draw[arrow] (repair) -- (select);
    \draw[arrow] (select) -- (check);

    \draw[arrow] (check.south) -- ++(0,-0.35) -- ++(-6.47, 0) |- (dnr.west);

    \node[draw=none, above right=-0.3cm and -2.5cm of ia_box] {Integrated approach$^*$};

    \node[draw=none, minimum width=0, minimum height=0] () at (-1.1, -7.6) {*Executed independently for each scenario $s \in S$.}; 
\end{tikzpicture}
    \caption{Schematic of the PH algorithm.}
    }
    \label{fig:ph_diagram}
\end{figure}
}

PH begins with a feasible incumbent solution that includes all candidate lines and updates it whenever a new feasible solution with an improved cost is found. At the end of each PH iteration, two solutions are generated. The $\beta$-intersection solution includes candidate circuits built in at least $\beta\%$ of the scenarios, whereas the union solution comprises all candidate lines selected in at least one scenario. The costs of both solutions are computed for each scenario, and if power flow violations occur, a repair method is applied (see Section \ref{sec:dnr}). Because repair can possibly fail to completely eliminate violations at a scenario, a penalized cost is adopted to compare feasible and infeasible solutions. Note that lines added to repair a specific scenario are also inserted for the remaining scenarios in the same solution to maintain consistency. PH then begins the next iteration with the \blue{most recent feasible solution}. The algorithm terminates when it reaches the time limit.

\subsection{Destroy-and-repair procedure}\label{sec:dnr}

D\&R starts with a feasible initial solution $I$, and iteratively attempts to remove candidate lines based on their residual flows, following a scheme similar to binary search. Algorithm \ref{alg:destroy_and_repair} sketches the D\&R procedure.

\setcounter{algocf}{0}
\renewcommand{\algorithmcfname}{Algorithm}
\LinesNumbered
\begin{algorithm}[htbp]
    \scriptsize
    \caption{\texttt{destroyAndRepair} procedure}
    \label{alg:destroy_and_repair}
    \algorithmfootnote{$^*$Computed as ${\scriptstyle (\bar{F}\,^1_k - |L_{f\,^1_k}|)/\bar{F}\,^1_k}$, where $L_{f^1_k}$ denote the value of variable \\ $f^1_k$ at the current solution.}
    \label{alg:dnr}
    \SetAlgoLined
\DontPrintSemicolon

\SetKwFunction{Select}{select}
\SetKwFunction{Sort}{sortByResidualFlows}
\SetKwFunction{Termination}{terminationCriteria}
\SetKwFunction{Min}{min}
\SetKwFunction{SolveLp}{solveLp}
\SetKwFunction{RemoveLines}{remove}
\SetKwFunction{AddLines}{insert}
\SetKwFunction{SolveModel}{solve}
\SetKwFunction{UpdateData}{updateSolution}
\SetKwFunction{Repair}{repair}    
\SetKwFunction{CompCost}{cost}
\SetKwFunction{CompViol}{violation}
\SetKwFunction{Time}{time}
\SetKwInOut{Input}{Input}
\SetKwInOut{Output}{Output}

\Input{Inserted \blue{$I$ and removed $R$} candidate lines; maximum number of iterations $it_{dr}$; time limit $tl$.}

\Output{Candidate lines to insert $I'$; candidate lines to remove $R'$.}

$rt \leftarrow 0.5, dt \leftarrow \blue{0.25}, rm \leftarrow \emptyset$\label{alg:dnr:initialization}

$\mathcal{K}'\leftarrow I$\label{alg:dnr:first_K_update}\;
$L \leftarrow$ \SolveLp{}\label{alg:dnr:init_update_model}

$c \leftarrow$  \CompCost{$L$}\label{alg:dnr:best}

$rm' \leftarrow$ \Select{$L, rt$}\tcp*{select by residual flows$^*$}\label{alg:dnr:init_rm_in}

\While{termination criteria is not reached}{\label{alg:dnr:beg_loop}
    $\mathcal{K}'\leftarrow\mathcal{K}'\setminus rm'$\label{alg:dnr:update_K_rm_lines}\;
    $L \leftarrow$ \SolveLp{}\label{alg:dnr:rm_lines}

    $v \leftarrow $\CompViol{$L$}\label{alg:dnr:comp_viol}

    \If{$v > 0$}{\label{alg:dnr:if_repair}
        $rm', v\leftarrow$\Repair{$rm', L, v$}\tcp*{Algorithm \ref{alg:repair}}\label{alg:dnr:repair}
    }
    $improved \leftarrow false$\;
    \If{$v = 0$}{\label{alg:dnr:if_feas}
        $c' \leftarrow$ \CompCost{$L$}

        \If{$c' < c$}{\label{alg:dnr:cost_impr}
            $improved \leftarrow true$
            
            $rm \leftarrow rm'$\label{alg:dnr:update_rm}
            
            $c \leftarrow c'$\label{alg:dnr:update_best}
            
            $rt \leftarrow rt + dt$\label{alg:dnr:increase_rt}
        }
    }

    \If{not $improved$}{\label{alg:dnr:beg_otherwise}
        $rt \leftarrow rt - dt$\label{alg:dnr:decrease_rt}\;
        $\mathcal{K}'\leftarrow I$\;
        $L \leftarrow$ \SolveLp{}\label{alg:dnr:reinsert_lines}\;
    }\label{alg:dnr:end_otherwise}

    $dt \leftarrow 0.5\times dt$
    
    $rm' \leftarrow$ \Select{$L, rt$}
}\label{alg:dnr:end_loop}

$I' \leftarrow I \setminus rm$\label{alg:dnr:updateI}

$R' \leftarrow R \cup rm$\label{alg:dnr:updateR}

\Return $I', R'$\label{alg:dnr:return}
\end{algorithm}

At the beginning of the algorithm, the LP model is updated according to parameter $I$, and solved (lines \ref{alg:dnr:first_K_update} and \ref{alg:dnr:init_update_model}, respectively). The best cost is then computed according to \eqref{lp_model:obj}, and including investment costs (line \ref{alg:dnr:best}). The set of lines to remove $rm'$ is created by selecting \blue{a ratio of $rt \in (0, 1]$} of the inserted candidate circuits with the largest residual flows (line \ref{alg:dnr:init_rm_in}). The main loop of \texttt{destroyAndRepair} (lines \ref{alg:dnr:beg_loop}--\ref{alg:dnr:end_loop}) is executed while the termination criteria is not met. We break the loop execution if its number of runs reaches the maximum number of iterations $it_{dr}$, the procedure's runtime reaches the time limit $tl$, $rm'$ is empty or its length between consecutive iterations remains unchanged. At each iteration, lines in $rm'$ are removed from the model (line \ref{alg:dnr:update_K_rm_lines}) and a new solution is generated by the solver (line \ref{alg:dnr:rm_lines}). \blue{The violation is computed as the sum of the values of the slack variables (line \ref{alg:dnr:comp_viol})}. In case the new restricted set $\mathcal{K}'$ generates a solution with power flow violations (line \ref{alg:dnr:if_repair}), circuits from $rm'$ can possibly be reinserted into $\mathcal{K}'$ by means of the \texttt{repair} operator (line \ref{alg:dnr:repair}). If $rm'$ yields a feasible solution (line \ref{alg:dnr:if_feas}) with an improved best cost (line \ref{alg:dnr:cost_impr}), the best $rm$ and $c$ are updated, and $rt$ increases (lines \ref{alg:dnr:update_rm}, \ref{alg:dnr:update_best}, and \ref{alg:dnr:increase_rt}, respectively). Otherwise, if $rm'$ results in an infeasible solution that cannot be fully repaired by the \texttt{repair} method, $rt$ is decreased and lines in $I$ are reinserted (lines \ref{alg:dnr:beg_otherwise}-\ref{alg:dnr:end_otherwise}). At the end of the loop, the new sets $I'$ and $R'$ are updated and returned (lines \ref{alg:dnr:updateI}--\ref{alg:dnr:return}).

As \texttt{destroyAndRepair} is based on the binary search algorithm, half the candidate lines with the largest residual flows are removed at the first iteration. If, for example, an improvement is found, $rt$ increases from $0.5$ to $0.75$; otherwise, $rt$ decreases to $0.25$. \blue{In both cases, $dt = 0.25$. If \texttt{destroyAndRepair} continues its execution, $rt$ is incremented or decremented by $dt = 0.125$ in the next iteration. The increment/decrement value $dt$ keeps reducing following the binary search scheme.}

The \texttt{repair} procedure is presented in algorithm \ref{alg:repair}.

\begin{algorithm}[htbp]
    \scriptsize
    \caption{\texttt{repair} procedure}
    \label{alg:repair}
    \SetAlgoLined
\DontPrintSemicolon

\SetKwFunction{Select}{selectWithViolation}
\SetKwFunction{AddLines}{solveLp}
\SetKwFunction{Round}{round}
\SetKwFunction{CompViol}{violation}
\SetKwFunction{Time}{time}
\SetKwInOut{Input}{Input}
\SetKwInOut{Output}{Output}

\Input{Removed candidate lines $rm$; current solution $L$; current violation $v$}

\Output{The best set of candidate lines to remove $rm$ and the associated violation $v$; }

\tcc{get lines in $rm$ whose associated existing line have flow violation}
$ri \leftarrow$\Select{$L, rm$}\label{alg:repair:init_ri}\;

\While{termination criteria is not reached}{\label{alg:repair:beg_loop}
    $\mathcal{K}'\leftarrow\mathcal{K}'\cup ri$\;
    $L \leftarrow$ \AddLines{}\label{alg:repair:add_lines}\;
    
    $v' \leftarrow $\CompViol{$L$}\label{alg:repair:comp_viol}\;
    \uIf{$v' < v$}{\label{alg:repair:if}
        $rm \leftarrow rm \setminus ri$\label{alg:repair:update_best_rm}\;
        $v \leftarrow v'$\label{alg:repair:update_best_viol}\;
        $ri \leftarrow$\Select{$L, rm$}\label{alg:repair:update_ri}\;
    }\Else{
        \textbf{break}\label{alg:repair:break}\;
    }
}\label{alg:repair:end_loop}

\Return $rm, v$\;\label{alg:repair:return}

\end{algorithm}

\texttt{repair} starts by selecting candidate lines in $rm$ whose corresponding existing line violates power flow constraints (line \ref{alg:repair:init_ri}). The main loop of \texttt{repair} (lines \ref{alg:repair:beg_loop}--\ref{alg:repair:end_loop}) stops when $v$ is equal to 0 or $ri$ is empty. At each iteration, a new solution is generated based on the candidate lines in $ri$ and the new violation is computed (lines \ref{alg:repair:add_lines} and \ref{alg:repair:comp_viol}, respectively). If the violation decreases (line \ref{alg:repair:if}), $rm$ and $v$ are updated (lines \ref{alg:repair:update_best_rm} and \ref{alg:repair:update_best_viol}, respectively), and new lines are selected to populate $ri$ (line \ref{alg:repair:update_ri}); otherwise, the loop stops (line \ref{alg:repair:break}). 
At the end, the algorithm returns the best set $rm$ and violation $v$ found (line \ref{alg:repair:return}).

\subsection{Beam search heuristic}

Our BS implementation is mostly based on \cite{velez2016beam, morais2024exact}, while incorporating several new features. In contrast to D\&R, BS \blue{discards in}feasible solutions \blue{without attempting to repair them}. The procedure begins with a feasible solution at the root node and iteratively explores a graph, one level at a time, where each path represents a possibly different solution. At each level, $N$ nodes are randomly selected for branching out of the $\lfloor(1 + \gamma)N\rfloor$ nodes with the best costs, computed according to \eqref{lp_model:obj} and including investment costs for newly added lines.

To create the branches, two strategies are applied sequentially. The first creates a random permutation of the inserted candidate lines to promote a search diversification. The second orders these circuits in non-increasing order of cost to intensify the search around promising regions. Each strategy partitions the inserted candidate circuits at the current node into subsets of $\max\{\eta{\scriptstyle\frac{|I|}{|\mathcal{K}|}\frac{|\mathcal{B}|}{1000}}, 1\}$ elements, where $\eta$ is an algorithm parameter, $|I|$ is the number of inserted lines at the current node, and $|\mathcal{K}|$ and $|\mathcal{B}|$ denote the total number of candidate circuits and system nodes, respectively. In the second strategy, lines whose removal fails to improve the parent node’s cost are excluded from future iterations. Each strategy continues until a fixed number $it_{bs}$ of consecutive levels are explored without improvement in the incumbent solution.

Each subset is assigned a cost $Z$, defined as the sum of the investment costs $C_k^{inv}$ of its candidate lines. For each node, up to $\omega$ branches are generated by removing the subsets of inserted circuits with the largest $Z$ values. If removing a subset produces a node with a larger objective value (or upper bound $UB$), the removed circuits are reinserted, and the node retains the same set of inserted candidates as its parent. However, the $UB$ of this new node is updated with the larger objective value, thereby imposing a penalty on it.

Fig. \ref{fig:bs_tree} illustrates a partial BS execution for an instance with 5 buses, 6 \blue{and 12} existing and candidate \blue{transmission lines, respectively}. For the sake of simplicity, let the sizes of subsets be set to 1, with parameters $\omega = 2$, $\gamma = 0.5$, and $N = 3$. This configuration implies that 3 out of the best $\lfloor(1 + 0.5)3\rfloor = 4$ subsets are selected for branching at each level.

\begin{figure*}[htbp]
    \centering
    \scalebox{0.70}[0.75]{ 
        \begin{tikzpicture}[edge from parent fork down,->,draw]

\begin{scope}[every node/.style={draw, rectangle, align=center}]
\node[level 0] (c0) {$\{1, 2, ..., 12\}$\\UB: $30521$}
    child {node [level 1] (c1) {Rm: $\{7\}$\\UB: $30419$}
        child {node [level 2] (c11) {Rm: $\{7,5\}$\\UB: $30318$}
            child {node [level 3] (c111) {Rm: $\{7,5\}$\\UB: $36569$}
            }
            child {node [level 3, pruned] (c112) {Rm: $\{7,5\}$\\UB: $36569$}
            }
        }
        child {node [level 2] (c12) {Rm: $\{7\}$\\UB: $32124$}
            child {node [level 2] (c121) {Rm: $\{7,5\}$\\UB: $30318$}
            }
            child {node [level 2] (c122) {Rm: $\{7\}$\\UB: $30320$}
            }
        }
    }
    child {node [level 1] (c2) {Rm: $\{5\}$\\UB: $30419$}
        child {node [level 2] (c21) {Rm: $\{5,7\}$\\UB: $30318$}
            child {node [level 3, pruned] (c211) {Rm: $\{5,7,4\}$\\UB: $30219$}
                }
            child {node [level 3, pruned, co] (c212) {Rm: $\{5,7\}$\\UB: $36569$}}
            }
        child {node [level 2, pruned] (c22) {Rm: $\{5,4\}$\\UB: $30320$}
            }
    };

    \node [draw=none,align=center](table) at (9,-4.95) {
        \scriptsize
        \setlength{\tabcolsep}{2pt} 
        \begin{tabular}{rrcrr}
            \toprule
            \multicolumn{1}{c}{\textbf{Line $\bm{k}$}} & \multicolumn{1}{c}{$\bm{C^{inv}_k}$} & & \multicolumn{1}{c}{\textbf{Line $\bm{k}$}} & \multicolumn{1}{c}{$\bm{C^{inv}_k}$} \\
            \cmidrule(){1-2} \cmidrule(){4-5}
            1  & 36  & & 7  & 101 \\
            2  & 99  & & 8  & 21  \\
            3  & 99  & & 9  & 21  \\
            4  & 99  & & 10 & 94  \\
            5  & 101 & & 11 & 94  \\
            6  & 99  & & 12 & 36 \\
            \bottomrule
        \end{tabular}
    };

    \path (c0) -- coordinate[midway] (c0c1) (c1);
    \draw node[above left, draw=none] at (c0c1) {\footnotesize $\{7\}$\\\footnotesize Z: $101$};
    \path (c0) -- coordinate[midway] (c0c2) (c2);
    \draw node[above right, draw=none] at (c0c2) {\footnotesize $\{5\}$\\\footnotesize Z: $101$};

    \path (c1) -- coordinate[midway] (c1c11) (c11);
    \draw node[above left, draw=none] at (c1c11) {\footnotesize $\{5\}$\\\footnotesize Z: $101$};
    \path (c1) -- coordinate[midway] (c1c12) (c12);
    \draw node[above right, draw=none] at (c1c12) {\footnotesize $\{2\}^*$\\\footnotesize Z: $99$};

    \path (c2) -- coordinate[midway] (c2c21) (c21);
    \draw node[above left, draw=none] at (c2c21) {\footnotesize $\{7\}$\\\footnotesize Z: $101$};
    \path (c2) -- coordinate[midway] (c2c22) (c22);
    \draw node[above right, draw=none] at (c2c22) {\footnotesize $\{4\}$\\\footnotesize Z: $99$};

    \path (c11) -- coordinate[midway] (c11c111) (c111);
    \draw node[above left, draw=none] at (c11c111) {\footnotesize $\{2\}^*$\\\footnotesize Z: $99$};
    \path (c11) -- coordinate[midway] (c11c112) (c112);
    \draw node[above right, draw=none] at (c11c112) {\footnotesize $\{3\}^*$\\\footnotesize Z: $99$};

    \path (c12) -- coordinate[midway] (c12c121) (c121);
    \draw node[above left, draw=none] at (c12c121) {\footnotesize $\{5\}$\\\footnotesize Z: $99$};
    \path (c12) -- coordinate[midway] (c12c122) (c122);
    \draw node[above right, draw=none] at (c12c122) {\footnotesize $\{4\}^*$\\\footnotesize Z: $99$};

    \path (c21) -- coordinate[midway] (c21c211) (c211);
    \draw node[above left, draw=none] at (c21c211) {\footnotesize $\{4\}$\\\footnotesize Z: $99$};
    \path (c21) -- coordinate[midway] (c21c212) (c212);
    \draw node[above right, draw=none] at (c21c212) {\footnotesize $\{3\}^*$\\\footnotesize Z: $99$};

    \draw node[draw=none] at (-1.15,0.5) {$\pi_1$};
    
    \draw node[draw=none] at (-7.1,-1) {$\pi_2$};
    \draw node[draw=none] at (4.95,-1) {$\pi_3$};
    
    \draw node[draw=none] at (-10.05,-2.5) {$\pi_4$};
    \draw node[draw=none] at (-4.05,-2.5) {$\pi_5$};
    \draw node[draw=none] at (1.95,-2.5) {$\pi_6$};
    \draw node[draw=none] at (7.95,-2.5) {$\pi_7$};

    \draw node[draw=none] at (-11.8,-4) {$\pi_8$};
    \draw node[draw=none] at (-8.3,-4) {$\pi_9$};
    \draw node[draw=none] at (-5.85,-4) {$\pi_{10}$};
    \draw node[draw=none] at (-2.35,-4) {$\pi_{11}$};
    \draw node[draw=none] at (0,-4) {$\pi_{12}$};
    \draw node[draw=none] at (3.65,-4) {$\pi_{13}$};

    \draw node[above left, draw=none, rotate=90] at (-10.58,-4.9) {\textbf{...}};
    \draw node[above left, draw=none, rotate=90] at (-4.58,-4.9) {\textbf{...}};
    \draw node[above left, draw=none, rotate=90] at (-1.08,-4.9) {\textbf{...}};

    \draw node[draw=none] at (-7.76,-5.6) {\footnotesize $^*$Indicates a setting that worsens the best objective found.};
\end{scope}
\end{tikzpicture}
    }
    \caption{Partial BS execution for an instance with 5 buses, 6 existing lines and 12 candidate circuits.}
    \label{fig:bs_tree}
\end{figure*}

The search begins at the root node $\pi_1$, which contains all candidate lines and an initial upper bound $UB = 30521$. Level 0 consists of 12 subsets, each containing a single circuit. After applying a random permutation, the $\omega = 2$ sets with the largest $Z$ values are selected (i.e., $\{7\}$ and $\{5\}$, both with $Z = 101$) which generate nodes $\pi_2$ and $\pi_3$ with the same $UB = 30419$.

At level 1, the remaining candidate circuits in each node are randomly permuted \blue{again} and divided into single-element subsets. For node $\pi_2$, the two subsets with the largest values of $Z$ are $\{5\}$ and $\{2\}$, with costs $101$ and $99$, respectively. On the one hand, removing subset $\{5\}$ improves the upper bound to $UB = 30318$. On the other hand, removing $\{2\}$ generates a new node with $UB = 32124$. In the latter, circuit 2 is reinserted into the solution. However, the larger upper bound is maintained. The procedure is applied similar\blue{ly} to node $\pi_3$, produc\blue{ing} two new nodes $\pi_6$ and $\pi_7$.

Level 2 contains 4 nodes: $\pi_4$, $\pi_5$, $\pi_6$, and $\pi_7$. From these, 3 out of the 4 nodes with the best corresponding upper bounds are randomly selected, namely, nodes $\pi_4$, $\pi_5$, and $\pi_6$. In this figure, a pruned node is indicated by two lines crossing it. This partial BS execution finishes with the best solution found \blue{at} node $\pi_{12}$, \blue{with} cost $30219$. Although this node is pruned in subsequent iterations, its cost and corresponding solution are stored for later \blue{use}. 

\section{Computational experiments}\label{sec:comp_experiments}

The proposed algorithms were implemented in Julia \citep{bezanson2017julia}, version 1.11.6, with the modeling language JuMP \citep{lubin2023jump}, whereas the mixed-integer and linear programs were modeled and solved with Gurobi 12.0.3. \blue{Except for the number of threads, all solver parameters were set to their default values.} All experiments were conducted \blue{on a server equipped with} an Intel\textsuperscript{\tiny\textregistered}~Xeon\textsuperscript{\tiny\textregistered}~Silver 4316 at 2.3 GHz, 128 GB of RAM\blue{, and 40 physical cores (80 hardware threads)}. Due to resource limitations, the nodes in the \blue{D\&R-}BS procedure were processed using a single thread. Moreover, we parallelized the PH subproblems with MPI. Our implementation is available at \url{https://github.com/yurerocha/TEP_PSCC2026.jl}\blue{, and the system data are available at \url{https://github.com/yurerocha/StochasticPGLibOPF.jl}.}

\blue{Our server provides 80 hardware threads. In our PH implementation, each subproblem is solved by a worker process, and each worker process uses a single thread. Solving all 96 subproblems simultaneously would therefore require 97 threads (96 for the MPI worker processes and one for the master process). Due to the limited number of hardware threads, we instead use 49 threads (48 worker processes and one master process), which results in approximately twice the runtime per PH iteration.}

\subsection{Benchmark instances}

To evaluate our methods, we adopted a few systems from the power grid library for OPF (PGLib-OPF) \citep{babaeinejadsarookolaee2019power} and the California test system \citep{taylor2024cats}. For the deterministic case, we adapted the PGLib-OPF systems as follows. For each transmission line in the original instance, we created two candidate lines with identical electrical parameters. Although identical, the two candidate lines represent different investment decisions. Existing lines with zero reactance were removed. Both demands and generation capacities (lower and upper bounds) were doubled to scale system size. \blue{The costs of candidate transmission lines were computed as described in \cite{musselman2025climate}, using cost estimates reported in \cite{MISO2024TransmissionCostGuide}. Generation costs were adjusted for inflation and expressed in 2024 USD. Moreover, random values from the set $\{0.00, 0.02, ..., 0.09, 0.10\}$ were assigned to break symmetry among candidate transmission lines corresponding to the same existing line.}

To construct multi-scenario instances, we extended the PGLib-OPF systems using parameters from the CATS \blue{dataset}. First, renewable generators were designated by iteratively selecting generators whose capacity deviated least from the average renewable generation capacity in the CATS base case. One generator was selected per iteration until the cumulative capacity of the newly assigned renewable generators, relative to \blue{the} total generation in the \blue{PGLib-OPF} system, reached the same proportion as in the CATS base case. Second, we generated multiple scenarios with varying renewable capacities and demands. Specifically, for each CATS scenario, generator capacities in the PGLib-OPF system were scaled according to the ratio of renewable generation with respect to the total generation in that scenario. All procedures involving generator capacities were applied separately for solar and wind renewable sources. Next, we scaled the demands so that total demand-to-generation ratio matches that of the corresponding CATS scenario.

Out of the $365\times24=8760$ hourly scenarios in CATS, the ones representing days March 20th, June 20th, September 22nd, and December 21st \blue{(representing one day for each season)} were extracted (\blue{yielding} a total of $4\times24=96$ hourly scenarios) each with probability $p_w$. \blue{Although investment decisions are static, these scenarios allow operational decisions to represent different operating conditions throughout the year. The criteria for selecting representative days are simple and practically motivated.} 

Candidate lines in the multi-scenario systems were generated as in the deterministic case. Since the modified scenarios sometimes resulted in total demand below twice \blue{the total demand} of the PGLib-OPF system, an additional scaling step was applied. \blue{The 96 hourly} scenarios were \blue{then} sorted in non-decreasing order of total demand, and the scaling factor required for the scenario at index $round(0.8\times|\mathcal{S}|)$ (i.e., $\approx$ \blue{the} 80th percentile) to reach twice the \blue{total PGLib-OPF system} demand was computed. This factor was then applied to scale generation capacities and demand \blue{levels} consistently across all scenarios. The lower bounds on generator capacities in the stochastic systems were set to zero, as in the CATS scenarios. Also, we performed preliminary tests to assert feasibility for solutions that build all candidate lines for the deterministic and the stochastic systems. 

For the computational experiments, we considered both deterministic and stochastic versions of systems \texttt{wpk3012}, \texttt{rte6495}, \texttt{epigrids7336}, \texttt{cats8870}, \texttt{goc9591}, and \texttt{goc10000}, adapted as previously described (the numbers in the names represent the number of buses in the instance).

\subsection{Parameter tuning}\label{sec:param_tuning}

Now we briefly explain the choices in the parameter tuning. LP penalty parameter $\lambda$ was computed for each instance as the maximum value of generation upper bound times its cost. This parameter requires special attention. Very large $\lambda$ values can lead to numerical instability, especially on larger systems, whereas very small values might yield the LP model to produce misleading results (for example, a solution with power flow violations that would otherwise be feasible if a larger $\lambda$ value was used). For the solution threshold parameter $\beta$ of PH, we conducted preliminary experiments with 0.25 and 0.5. The latter means that at least half the scenarios agreed on the construction of a candidate line. This value seemed, however, too aggressive, constantly leading to infeasible solutions. Therefore, we adopted $\beta = 0.25$. 

The value selected for the parameter $it_{dr}$ of D\&R was 15, after experimenting with values 5, 10, and 15. In fact, $it_{dr} = 10$ showed a better compromise between computational time and solution quality for the deterministic cases, but $it_{dr} = 15$ had a positive impact on larger stochastic instances. BS parameters $\omega$, $N$, and $\gamma$, were the same as in \citet{morais2024exact}, whereas $it_{bs} = 15$. Larger values for this last parameter can result in too many BS \blue{levels} without improvements. Finally, we performed experiments with different values for $\eta$, namely, 0.0025, 0.005, and 0.01. The most promising results were obtained with $\eta = 0.005$.

\subsection{Integrated approach performance for the deterministic systems}

To better evaluate the performance of the heuristic approach, we conducted a set of experiments with the deterministic systems under a 10-minute time limit. Two settings were then adopted: MIP-TEP, which is the MIP solver with a warm-start solution that inserts all candidate lines; and D\&R-BS+MIP-TEP, which first runs the D\&R-BS heuristic and then executes the commercial solver. \blue{Specifically, in the latter, D\&R-BS is run for up to 10 minutes and the solver is executed for the remaining time, that is, 10 minutes minus the runtime of the D\&R-BS heuristic.}

Table \ref{tab:results_det} summarizes the results obtained. In this table, \blue{``Start'',} ``UB'' and ``D\&R-BS (s)'' refer to the \blue{cost of the warm-start ``build-everything'' solution}, the best upper bound found, and the time taken by the heuristic procedures, respectively. Furthermore, \blue{``Gap (\%)'' highlights the optimality gap achieved by the solver, which considers both the best upper and lower bounds found, while ``Imp (\%)'' reports the relative improvement upon the upper bound and is computed as follows.}

{
    \footnotesize
    \blue{
    \begin{align*}
        \text{Imp (\%)} = 100\times \frac{\text{UB}\{\text{MIP-TEP}\}-\text{UB}\{\text{D\&R-BS+MIP-TEP}\}}{\text{UB}\{\text{MIP-TEP}\}}
    \end{align*}
    }
}
\begin{table}[htbp]
    \centering
    \scriptsize
    \caption{Impact of the integrated approach for deterministic baseline scenarios under a 10-minute time limit.}
    \setlength{\tabcolsep}{2.0pt} 
    \resizebox{\linewidth}{!}{
    \begin{tabular}{lrrrrrrr}
        \toprule
        \multirow{2}{*}{\textbf{System}} & \multicolumn{1}{c}{\multirow{2}{*}{\blue{Start}}} & \multicolumn{2}{c}{\textbf{MIP-TEP}} & \multicolumn{4}{c}{\textbf{D\&R-BS+MIP-TEP}} \\
        \cmidrule(lr){3-4}\cmidrule(lr){5-8}
        & & \multicolumn{1}{c}{UB} & \multicolumn{1}{c}{\blue{Gap (\%)}} & \multicolumn{1}{c}{\blue{D\&R-BS (s)}} & \multicolumn{1}{c}{UB} & \multicolumn{1}{c}{\blue{Gap (\%)}} & \multicolumn{1}{c}{\blue{Imp (\%)}} \\
        \midrule
        \texttt{wpk3012} & \blue{7923243} & \blue{7656580} & \blue{2.74} & \blue{21.51} & \blue{7518058} & \blue{0.95} & \blue{1.81} \\
        \texttt{rte6495} & \blue{9521005} & \blue{9409689} & \blue{26.98} & \blue{61.63} & \blue{7188622} & \blue{4.42} & \blue{23.60} \\
        \texttt{epigrids7336} & \blue{5110551} & \blue{5083300} & \blue{32.43} & \blue{204.80} & \blue{4549212} & \blue{24.50} & \blue{10.51} \\
        \texttt{cats8870} & \blue{1723141} & \blue{1639536} & \blue{23.82} & \blue{197.73} & \blue{1345788} & \blue{7.19} & \blue{17.92} \\
        \texttt{goc9591} & \blue{3799469} & \blue{3757820} & \blue{64.86} & \blue{91.70} & \blue{2453576} & \blue{46.17} & \blue{34.71} \\
        \texttt{goc10000} & \blue{4604347} & \blue{4561517} & \blue{30.47} & \blue{128.97} & \blue{3344869} & \blue{5.19} & \blue{26.67} \\
        \midrule
        \textbf{Average} & \blue{5446959} & \blue{5351407} & \blue{30.22} & \blue{117.72} & \blue{4400021} & \blue{14.74} & \blue{19.20} \\
        \bottomrule
    \end{tabular}
    }
    \label{tab:results_det}
\end{table}
\blue{The} D\&R-BS+MIP-TEP \blue{framework} outperformed the MIP solver with a warm start for all systems. \blue{The framework achieved superior results, with an average cost improvement of 19.20\% and a reduction in the average optimality gap from 30.22\% to 14.74\%. In particular, the maximum time spent by the heuristic procedures was 204.80 seconds for system \texttt{epigrids7336}, whereas the largest gap reduction was achieved for system \texttt{goc10000}. This suggests that heuristic runtimes depend not only on the number of buses but also on the system structure. Furthermore, improvements tend to be more significant for larger systems.}

\blue{We also performed experiments without warm starts. In this case, the MIP solver failed to find an incumbent solution for any system within the configured time limit. As shown in Table \ref{tab:results_det}, when warm-start solutions were provided, Gurobi managed to improve them. In fact, the solver reduced the average cost to $5351407$, which corresponds to approximately $98.25\%$ of the average initial cutoff value.}

\subsection{\blue{Performance of the solution approach for the stochastic systems}}

To evaluate the performance of PH \blue{combined with the IA heuristic (PH$_{\text{IA}}$)}, we implemented a baseline algorithm that applies the integrated approach \blue{(BA$_{\text{IA}}$)} once for each subproblem on a single thread. The BA\blue{$_{\text{IA}}$ framework} starts from \blue{the} feasible solution \blue{in which} all candidate lines \blue{are constructed} and, at the end of its execution, selects the best \blue{of} two solutions, i.e., $\beta$-intersection and union, created in a similar fashion as described for the PH procedure. \blue {The} PH\blue{$_{\text{IA}}$} \blue{approach} executes multiple iterations with a time limit ($tl$) of 5 minutes for systems with up to 9000 buses and 10 minutes for larger systems. The overall time limit for PH\blue{$_{\text{IA}}$} is 2 hours. For BA\blue{$_{\text{IA}}$}, the time limit is set to 1 hour for all systems. Due to the relation between \blue{hardware} threads \blue{and} scenarios (see Section \ref{sec:param_tuning}), BA\blue{$_{\text{IA}}$} also runs for 2 hours. \blue{Furthermore, we also implemented an alternative version of the PH algorithm that employs the MIP solver on the subproblems (PH$_{\text{MIP-TEP}}$) and adopts the same parameters as PH$_{\text{IA}}$.} The results are summarized in Table \ref{tab:results_sto}. \blue{Here, ``Imp (\%)'' is computed as in the deterministic setting, but considering the procedures BA$_{\text{IA}}$ and PH$_{\text{IA}}$.}

\begin{table}[htbp]
    \centering
    \scriptsize
    \caption{Results obtained for the stochastic systems with a 2-hour time limit.}
    \begin{tabular}{lcccrr}
        \toprule
        \multirow{2}{*}{\textbf{System}} & \multirow{2}{*}{\textbf{Start}} & \multirow{2}{*}{\textbf{PH$_{\text{MIP-TEP}}$}} & \multirow{2}{*}{\textbf{BA$_{\text{IA}}$}} & \multicolumn{2}{c}{\textbf{PH$_{\text{IA}}$}} \\
        \cmidrule(lr){5-6}
        & & & & \multicolumn{1}{c}{UB} & \multicolumn{1}{c}{Imp (\%)} \\
        \midrule
        \texttt{wpk3012}        & 5709815 & 5681549 & 5603122 & 5505519 & 1.74 \\
        \texttt{rte6495}        & 6200923 & 6096518 & 4611551 & 4210422 & 8.70 \\
        \texttt{epigrids7336}   & 2844713 & 2817832 & 2531838 & 2363259 & 6.66 \\
        \texttt{cats8870}       & 1437393 & 1356825 & 1109278 & 995340 & 10.27 \\
        \texttt{goc9591}        & 1966510 & 1924856 & 759335 & 740091 & 2.53 \\
        \texttt{goc10000}       & 2746537 & 2701498 & 1490285 & 1463878 & 1.77 \\
        \midrule
        \textbf{Average}        & 3484315 & 3429846 & 2684235 & 2546418 & 5.28 \\
        \bottomrule
    \end{tabular}
    \label{tab:results_sto}
\end{table}

\blue{Table \ref{tab:results_sto} shows that, similarly to the deterministic systems, the MIP solver (here embedded within the PH algorithm) enhanced warm-start solutions in the PH$_{\text{MIP-TEP}}$ setting. However, the average cost decrease was minor compared to that achieved by PH$_{\text{IA}}$. Specifically, PH$_{\text{IA}}$ improved the incumbent solution in 60.66\% of its iterations on average. In contrast, PH$_{\text{MIP-TEP}}$ succeeded in 15.40\% of its iterations, which further demonstrates the impact of the integrated approach. In fact, incorporating our IA into the BA significantly enhanced the initial solutions compared to the more elaborate PH$_{\text{MIP-TEP}}$ approach. Moreover, the hybrid PH$_{\text{IA}}$ framework further increased robustness and effectiveness in producing high-quality solutions compared to BA$_{\text{IA}}$. In this case, PH$_{\text{IA}}$ achieved a 5.28\% improvement. These results confirm the effectiveness of both the IA heuristic and the PH algorithm for the stochastic instances.}

\section{Concluding remarks and future work}\label{sec:conclusion}

This paper propose\blue{s} a hybrid approach, based on D\&R, BS, and PH, for the TEP problem under uncertaint\blue{y}. Each subproblem \blue{is} solved by means of an integrated method that sequentially employs two heuristic procedures and the Gurobi solver. The first heuristic, D\&R, starts from an initial solution with constructed candidate lines. \blue{Next, it} quickly obtains feasible solutions by \blue{iteratively} removing candidate lines based on their residual power flows, eventually repairing solutions with power flow violations. The second heuristic, BS, simultaneously investigates multiple candidate solutions by attempting to remove various sets of candidate lines associated with large costs. Finally, the PH algorithm \blue{decomposes} large-scale stochastic problems into single-scenario subproblems, \blue{which are} linearized by means of the DC-OPF approximation. 

Computational experiments on adapted versions of \blue{the} systems \texttt{wpk3012}, \texttt{rte6495}, \texttt{epigrids7336}, \texttt{cats8870}, \texttt{goc9591}, and \texttt{goc10000} \blue{demonstrate} the algorithm's performance on both deterministic and stochastic \blue{settings}. In the \blue{deterministic case, the integrated approach, combining the} D\&R-BS \blue{operator with} an exact \blue{solver} \blue{outperforms} a \blue{stand-alone} commercial solver under a 10-minute time limit, \blue{achieving} an average \blue{cost improvement} of $\blue{19.20}\%$. \blue{Moreover}, the optimality gaps reported by Gurobi decreased from $\blue{30.22}\%$ to $\blue{14.74}\%$. In the \blue{stochastic case}, the PH\blue{$_{\text{IA}}$ framework} outperforms \blue{both BA$_{\text{IA}}$ and PH$_{\text{MIP-TEP}}$, achieving} an average \blue{cost} improvement of \blue{5.28\% over BA$_{\text{IA}}$.}

Prospective paths of research include the integration of regularization techniques such as Benders decomposition into the hybrid framework, as well as the extension to multi-stage stochastic optimization. \blue{Additional research directions include scenario clustering techniques for selecting and grouping representative days, incorporating uncertainty in generation expansion planning, and integrating energy storage systems into the model.}



%

{
    \footnotesize
    \bibliography{references}
}

\end{document}